\documentclass[aps,pre,groupedaddress,twocolumn]{revtex4-2}

\usepackage{graphicx}
\usepackage{multirow}
\usepackage{amsmath}
\usepackage{amssymb}
\usepackage{enumitem}
\usepackage{bm}
\usepackage[ruled,vlined]{algorithm2e}
\usepackage{hyperref}

\newcommand{\Cov}{\mathrm{Cov}}
\newcommand{\dist}{\mathrm{dist}}

\begin{document}


\title{Inferring Causal Networks of Dynamical Systems \\ through Transient Dynamics and Perturbation}


\author{George Stepaniants}
\email[]{gstepan@mit.edu}
\affiliation{Department of Mathematics, Massachusetts Institute of Technology, Cambridge, MA\\
Department of Mathematics, University of Washington, Seattle, WA}

\author{Bingni W. Brunton}
\email[]{bbrunton@uw.edu}
\affiliation{Department of Biology, University of Washington, Seattle, WA}

\author{J. Nathan Kutz}
\email[]{kutz@uw.edu}
\affiliation{Department of Applied Mathematics, University of Washington, Seattle, WA}



\begin{abstract}
Inferring causal relations from time series measurements is an ill-posed mathematical problem, where typically an infinite number of potential solutions can reproduce the given data. We explore in depth a strategy to disambiguate between possible underlying causal networks by perturbing the network, where the actuations are either targeted or applied at random. The resulting transient dynamics provide the critical information necessary to infer causality. Two methods are shown to provide accurate causal reconstructions:  Granger causality (GC) with perturbations, and our proposed {\em perturbation cascade inference} (PCI). Perturbed GC is capable of inferring smaller networks under low coupling strength regimes. Our proposed PCI method demonstrated consistently strong performance in inferring causal relations for small (2--5 node) and large (10--20 node) networks,  with both linear and nonlinear dynamics. Thus the ability to apply a large and diverse set of perturbations/actuations to the network is critical for successfully and accurately determining causal relations and disambiguating between various viable networks.
\end{abstract}


\maketitle

\section{Introduction}

The ability to determine causal relationships in complex, dynamical networks from time-series measurements alone is an important open challenge in the engineering, biological and physical sciences.  This task is challenging because inferring causal networks from time series observations alone is an ill-posed mathematical problem, and a potentially infinitude of solutions may accurately reproduce the given data.  Despite decades of research effort, the large and diverse set of mathematical methods that have been developed still have limitations in accurately reproducing causal network structures, especially for nonlinear dynamical systems on networks~\cite{lusch2016, sugihara2012,monster2017,daniella,emily,krakovska,angulo2017}.    

In certain applications, one has the ability to actuate the dynamical network of interest and generate additional information about the unknown network structure.  In this work, we exploit this capacity to perturb a dynamical network, in targeted or random ways, to extract an accurate directed graph of the true underlying dynamical system.  We demonstrate the accuracy and efficiency of our \emph{perturbation cascade inference} (PCI) algorithm on a variety of test problems, showing in which cases such network perturbation strategies can resolve this ill-posed inference problem.

Wiener first proposed a statistical notion of causality~\cite{wiener1956} by noting that $Y$ causes $X$ if knowing the past of $Y$ improves the prediction of $X$,  compared to using the past of $X$ alone.   Many innovations, modifications and reformulations of the causal inference problem have since been proposed~\cite{imbensrubin2015, pearl2009_1, morganwinship2014, holland1986}, with various statistical regularizations used to achieve potential solutions of this ill-posed problem.  The Nobel prize winning work of Granger~\cite{granger1969, granger1980} built upon Wiener's definition and formalized time series causality through linear regression of stochastic processes.  However, determining causal relationships in large networks of nonlinear dynamical systems still remains problematic with not only {\em Granger Causality} (GC)~\cite{lusch2016, sugihara2012}, but also with {\em convergent cross mapping} (CCM)~\cite{monster2017}, and other sparsity-promoting techniques (e.g. Ref.~\cite{daniella} Fig. 1, 2 and Ref.~\cite{emily} Fig. 4). Krakovsk\'a et al~\cite{krakovska} provide useful inference methods on a very limited set of systems and fail to produce even approximately accurate directed graphs for more complex network dynamics.  

\begin{figure}[t]
\includegraphics[width=9cm]{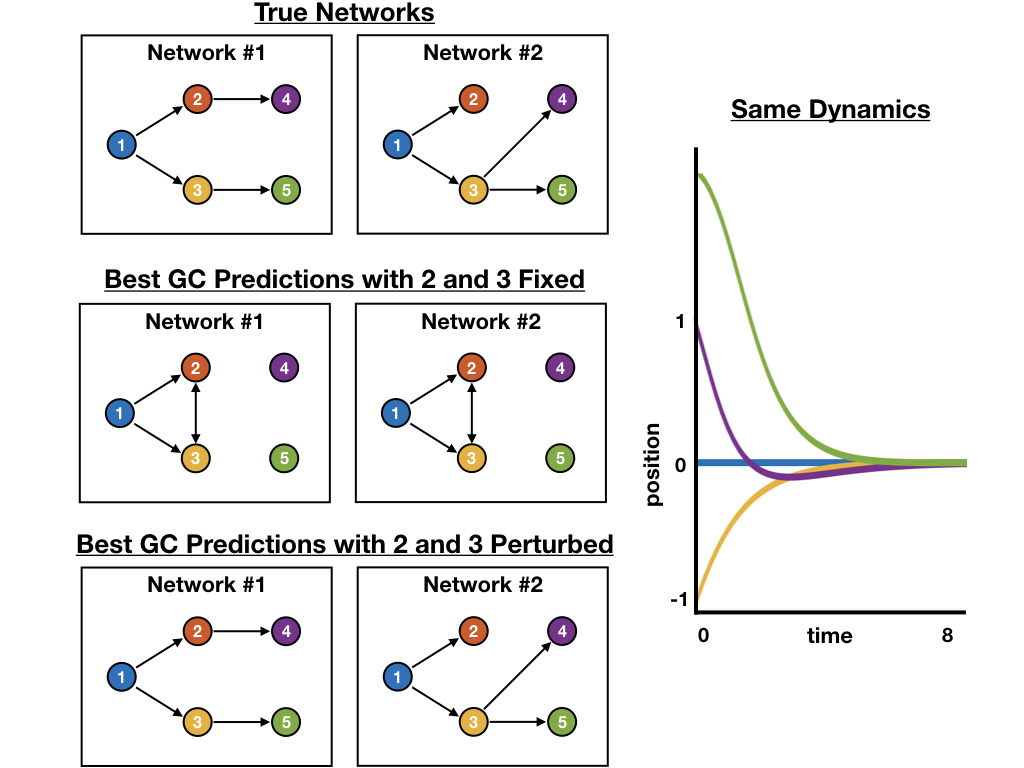}%
\caption{One simple example to demonstrate the fundamentally ill-posed nature of network inference. Specifically, two Kuramoto networks can easily produce identical time-series data, and it is impossible to  disambiguate between them from this data alone.
We show the best-case predictions of GC. When nodes 2 and 3 are perturbed (i.e. their initial conditions are allowed to vary), GC is able to perfectly resolve both graphs.}
\label{fig:netinfProblems}
\end{figure}

\begin{figure*}[t]
\includegraphics[width=15cm]{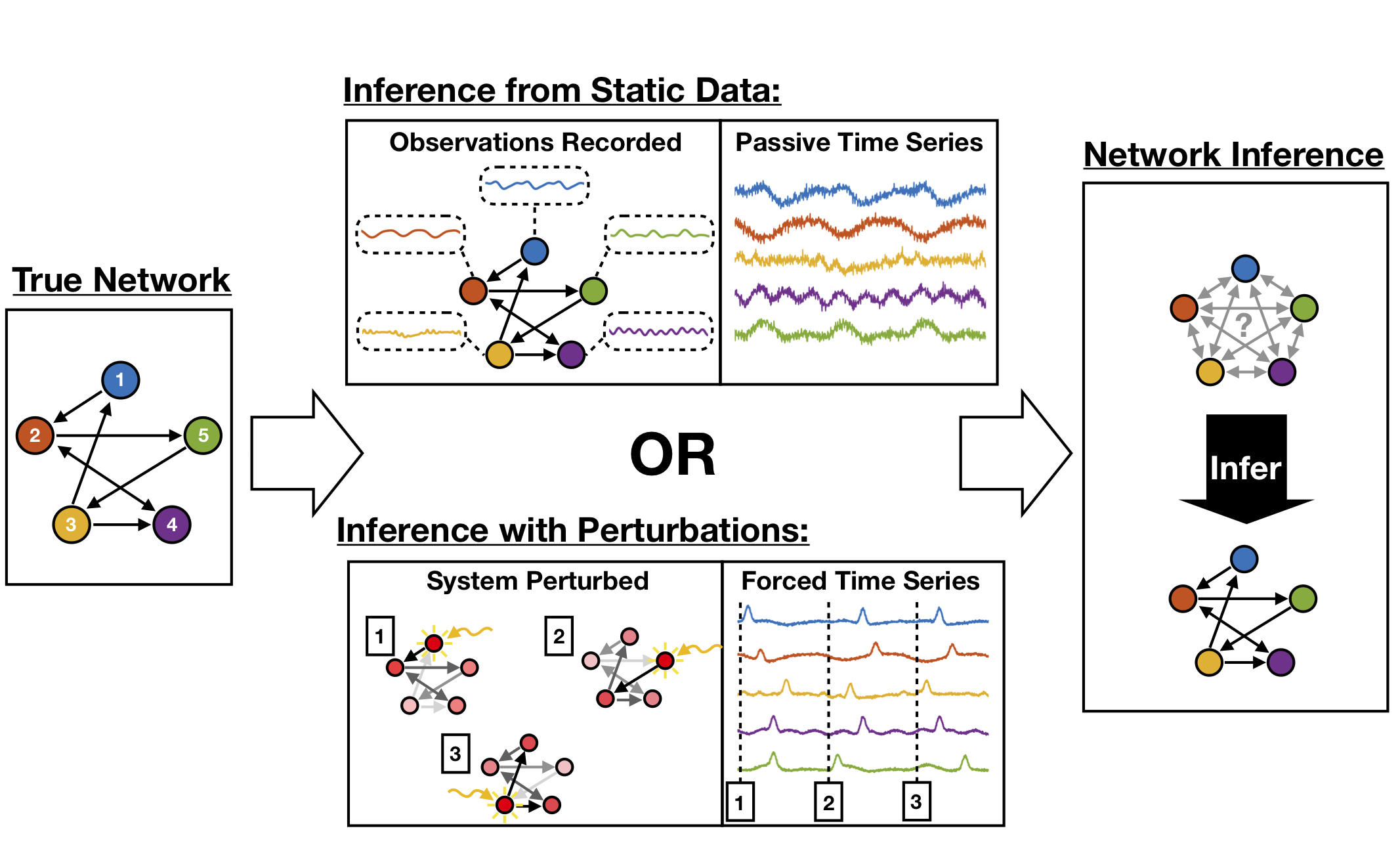}%
\caption{An illustration of two approaches to network inference. Given the true network shown on the left, one approach is to passively collect time series observations of each network node (top center) and use time series analysis to infer the network structure shown on the right. If perturbations of the network nodes are allowed, another approach is to systematically force the network nodes and study their dynamics after each perturbation to infer the network structure.}
\label{fig:inferenceSummary}
\end{figure*}

Fig.~\ref{fig:netinfProblems} illustrates how inferring networks from observational data is an ill-posed problem. In particular, identical time series data can be collected from several different network topologies, making disambiguation of underlying network structures impossible on the given data.  In such cases, it is necessary to push the networked dynamical system into regimes where its dynamics are no longer degenerate; in other words, we need observations in regimes where it is not possible for two different networks to produce the same observed trajectories.

Leveraging the ability to actuate or perturb components of a networked system can lead to significant insights about its network topology and dynamics~\cite{timme2014}.  Perturbation inference methods subject a networked dynamical system to external driving forces and measure its collective response to reverse-engineer the underlying structure.  For instance, drug perturbations and gene knockouts are applied extensively in the study of gene regulatory networks~\cite{ud2016optimal, meinshausen2016methods}, and activations of genetically modified neurons in optogenetics enable the identification of specific neuronal connectivities~\cite{lepperod2018inferring}. A large body of research has also investigated inference of networked systems subject to small perturbations about a stable point, particularly in the context of gene regulatory network reconstruction~\cite{bansal2006, bansal2007, gardner2003}.  By observing response dynamics of a system from a series of systematic perturbations, various inference methods have been able to model the time dynamics and underlying network topology of coupled oscillator systems, cellular signaling networks and competitive economic markets~\cite{timme2007, molinelli2013, kuypers2012, delabays2020network}.

The present work discusses approaches to reverse engineer the network connectivity  by perturbing, or actuating, the networked dynamical system as shown in Fig.~\ref{fig:inferenceSummary}.  We use simulated nonlinear dynamical systems, for which the ground-truth network structure is known, to study the network reconstruction quality of the various inference methods in its passive, unperturbed state, and under the influence of actuation or perturbation.  We propose two methods for uniquely disambiguating data for inference:  (a) judiciously perturbed GC with sampling only of the transient dynamics, and (ii) an active inference approach called \textit{Perturbation Cascade Inference} (PCI).
PCI infers the structure of a dynamical network through systematic perturbations of its components.  This approach is similar in spirit to {\em perturb and observe} methods for system identification for control~\cite{billings1980identification}.   In order to evaluate the performance of network inference algorithms in both the passive and active regimes, we simulate respectively the unperturbed and perturbed time dynamics of coupled spring-mass and Kuramoto oscillators on random Erdos-Renyi networks with varying edge densities, coupling strengths and sizes.  We show that by either targeted perturbations, or random global network actuation, we can accurately reconstruct the causal network with a sufficiently large number of perturbations.

\section{Background}
Before developing the proposed PCI algorithm, we briefly review the mathematical architecture of one of the first time-series causal inference techniques.  
This section also serves to establish some basic notions of causal relationships. 
Granger Causality (GC) infers the causal relationship between two multivariate time series $\mathbf{X}_t, \mathbf{Y}_t$ by fitting a vector autoregressive (VAR) model to test if $\mathbf{X}_t$ conditional on its own past does not depend on the past of $\mathbf{Y}_t$. Let $\mathbf{X}_t, \mathbf{Y}_t$ be $n \times 1$ vectors. A $p$th order VAR model for $\mathbf{X}_t$ has the form
\begin{equation}\label{eq:1}
\mathbf{X}_t = \sum_{k=1}^p \mathbf{A}_k\mathbf{X}_{t - k} + \bm{\epsilon}_t,
\end{equation}
where the $n \times n$ matrices $\mathbf{A}_k$ are the regression coefficients that specify how $\mathbf{X}_t$ depends on its past and the $n \times 1$ vector $\bm{\epsilon}_t$ contains the residuals or prediction error of the model.  GC tests whether past values of $\mathbf{Y}_t$ improve the prediction of $\mathbf{X}_t$ as compared to the past values of $\mathbf{X}_t$ alone. The joint regression of $\mathbf{X}_t$ on its own past and the past of $\mathbf{Y}_t$ is given by
\begin{equation}\label{eq:2}
\mathbf{X}_t = \sum_{k=1}^p \mathbf{A}'_k\mathbf{X}_{t - k} + \sum_{k=1}^p \mathbf{B}'_k\mathbf{Y}_{t - k} + \bm{\epsilon}'_t,
\end{equation}
where $\mathbf{A}'_k, \mathbf{B}'_k$ are $n \times n$ matrices and the residual $\bm{\epsilon}'_t$ is an $n \times 1$ vector. Using the residuals from both VAR models in (\ref{eq:1}) and (\ref{eq:2}), we can construct the log-likelihood ratio
$\mathcal{F}_{\mathbf{Y} \to \mathbf{X}} = \ln ({|\Sigma'|}/{|\Sigma|})$
where $\Sigma' = \Cov(\bm{\epsilon}'_t)$ and $\Sigma = \Cov(\bm{\epsilon}_t)$ are the covariance matrices of the residuals of both VAR models. This log-likelihood ratio is called the \textit{G-causality} from $\mathbf{Y}$ to $\mathbf{X}$ and it characterizes the causal influence of $\mathbf{Y}$ on $\mathbf{X}$. $\mathcal{F}_{\mathbf{Y} \to \mathbf{X}}$ is the test statistic for the null hypothesis of zero causality
\begin{equation}\label{eq:4}
\mathbf{B}_1' = \mathbf{B}_2' = ... = \mathbf{B}_p' = 0.
\end{equation}
We say that $\mathbf{Y}$ \textit{Granger causes} $\mathbf{X}$ if the G-causality $F_{\mathbf{Y} \to \mathbf{X}}$ is statistically significant.

In unconditional G-causality described above, we can erroneously infer that $\mathbf{Y}$ causes $\mathbf{X}$ even if there is no causal relationship, since both variables are dependent on a third latent, confounding variable $\mathbf{Z}$. The method of \textit{conditional G-causality} eliminates such spurious causalities by ``conditioning out" all potential common dependencies. In the context of network inference, we can apply conditional G-causality to determine the existence of a causal relationship between every pair of nodes in a network. This approach is known as \textit{pairwise-conditional G-causality} and is the basis of the Multivariate Granger Causality (MVGC) toolbox developed by Barnett \& Seth \textit{et al.} which we use in this paper \cite{barnettseth2014}.

\begin{figure}[t]
\includegraphics[width=9cm]{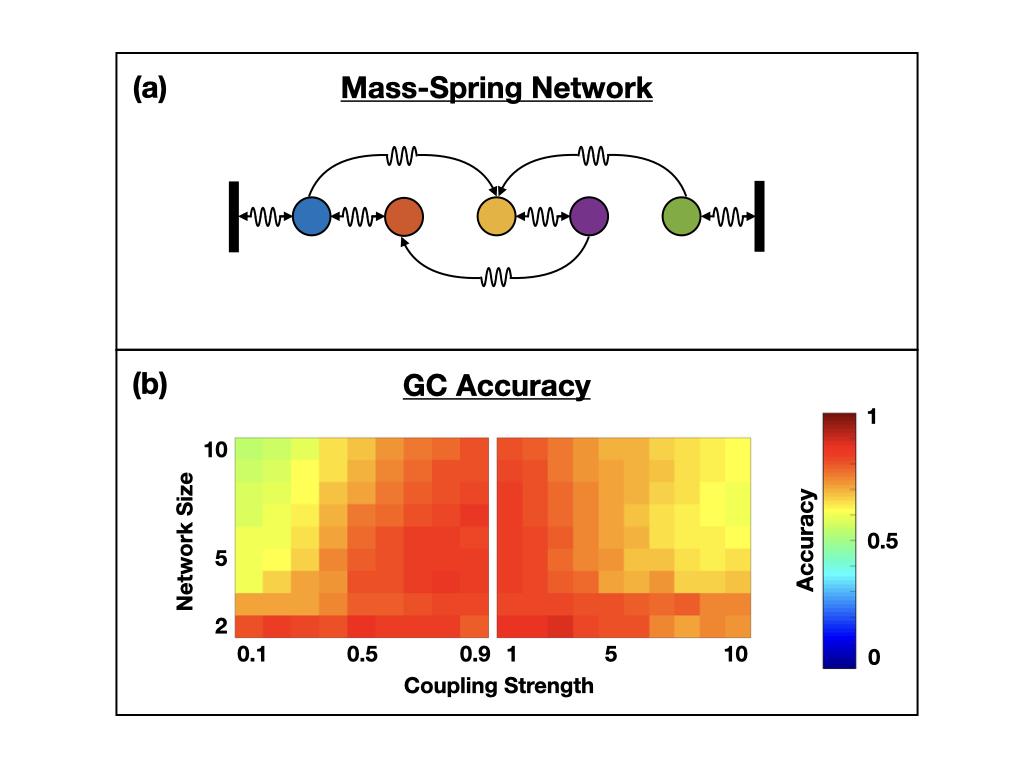}%
\caption{\textbf{(a)} An example of a 5 node mass-spring network with the endpoints attached to two fixed walls. \textbf{(b)} The performance of GC on such mass-spring networks varies with network sizes and coupling strengths. For coupling strengths near 1, GC is capable of inferring the network structure of the mass-spring system with 80-85\% accuracy for network sizes between 2 to 10 nodes. GC is only able to achieve accuracy above 85\% when the network size is 2. Furthermore, if we increase or decrease the coupling strength of the system away from 1, the inference accuracy of GC drops below 80\% for mass-spring networks of all sizes.}
\label{fig:GCHarmonicUnforced}
\end{figure}

\begin{figure*}[t]
    \centering
    \vspace*{-.5in}
    \includegraphics[width=15cm]{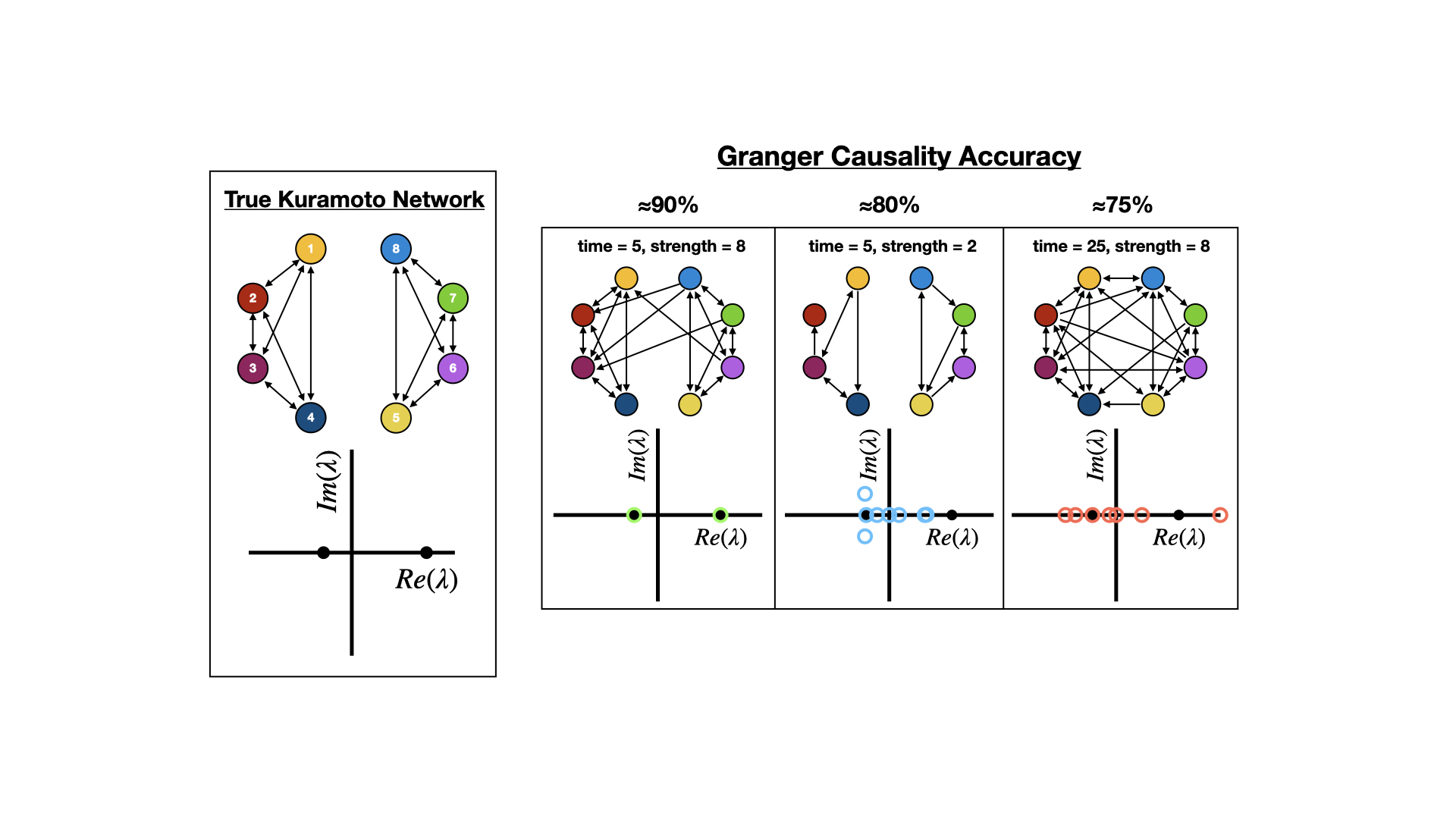}
    \vspace*{-.6in}
    \caption{Comparison of Granger inferred networks to the true underlying Kuramoto oscillator network with varied simulation endtimes and coupling strengths. Seemingly large edge prediction accuracies do not resemble the true underlying network structure. Moreover, the eigenvalues of the network connectivity matrix are not well characterized. Therefore, it would not be possible to use an inferred network that is 75-80\% accurate to simulate the time dynamics of a system.}
    \label{fig:inferenceAccuracy}
\end{figure*}

Results using GC inferences from the MVGC toolbox on time series observations of unforced mass-spring networks are shown in Fig.~\ref{fig:GCHarmonicUnforced}. Since GC performs an autoregressive linear fit to the observed time series data, it is expected to successfully infer the causal relationships in a linear dynamical system. Accordingly, GC makes highly accurate inferences on the mass-spring system for a wide range of network sizes when the coupling strengths (spring constants) are close to 1. The performance of this method rapidly degrades if the coupling strength is significantly smaller or larger than this value.  However, we are primarily interested in networked nonlinear dynamics, as most complex systems of interest are rarely linear.  Figure~\ref{fig:inferenceAccuracy} considers a network of nonlinear oscillators (Kuramoto oscillator network) and shows that MVGC fails to correctly identify the causal structure under most conditions~\cite{lusch2016}. This example motivates our need to develop a more robust method to extract causal relations.

There are a number of alternative approaches to GC for computing causality.  For instance, Convergent Cross Mapping ( CCM,~\cite{sugihara2012} assumes that the time-series data is sampled from a dynamical attractor and uses time-delay embeddings to construct attractor manifolds from the time-series observations of variables ${\bf X}$ and ${\bf Y}$. The correlation between the true state of ${\bf Y}$ and the predictions of ${\bf Y}$ given the ${\bf X}$-manifold determine the strength of the causal relationship from ${\bf Y}$ to ${\bf X}$.  Using information-theoretical approaches\cite{rahimzamani, krakovska}, it is possible to reconstruct causal graphs by computing statistics such as transfer entropy or directed information between the time series measurements of every pair of variables in a system. Similar to CCM, such approaches make no assumptions about the data generation process and often utilize time-lagged representations in the analysis.  Finally, a variety of model-based approaches have been proposed to fit the time series data through time-lagged regressions and ordinary differential equations. As described in \cite{daniella}, such methods often promote sparsity in the network parameters in order to restrict the space of possible solutions (e.g. obtain networks with the fewest number of edges that closely predict the system dynamics). Yet other methodologies such as Dynamic Causal Modeling \cite{friston} and Bayesian networks \cite{yu2004} have been developed to solve this ill-posed problem in application to neuronal dynamics and gene regulatory networks. Each method imposes a different regularization to extract a solution to the causal inference problem. However, these methods still face the same problem for causal inference in networked, nonlinear dynamics problems as illustrated in Fig.~\ref{fig:inferenceAccuracy}. In this paper, we show how the systematic perturbation of nodes in the network, either targeted or random, allow the resulting transient time dynamics to disambiguate the causal relations.

\section{Transients for Inference}
We will show how perturbations and transients can be used to reverse-engineer the structure of a time-dependent network.  Importantly, a distinct and diverse set of perturbations such as individual actuations of nodes or random kicks to the system are typically required in order to fully disambiguate the causal pathways.  Suppose we have a mass-spring network of five nodes with \textit{directed interactions} like the one depicted in Fig. \ref{fig:GCHarmonicUnforced}(a).
To learn about the connectivity of this system, we could perturb the blue node in the network away from its equilibrium position and observe how all the other oscillators respond to this perturbation over time.  After the initial mass is released a longitudinal wave travels through the network, displacing the masses and allowing us to observe the order in which they become displaced. Under the \textit{strict assumption} that all the coupling strengths between coupled oscillators are equal, we in fact observe that the red and yellow nodes become perturbed directly after the blue node, which subsequently activates the purple node. Note that there is no directed path from the blue to the green node, hence the green node remains unperturbed. This highlights, at least in part, the features exploited by algorithmic structures to infer causal relations.

\subsection{GC and Transients} 

Our first use of perturbations of a networked dynamical systems is with the standard GC framework for inference.  But instead of simply giving GC time-series measurements, we perturb the system by giving it different starting \textit{initial conditions}, and use the transient information to disambiguate the structured network.  As has already been shown, GC by itself fails to produce meaningful results in nonlinear systems unless time series lengths and coupling strengths of the system are chosen judiciously (see Fig.~\ref{fig:inferenceAccuracy}).

\begin{figure}
\includegraphics[width=8cm]{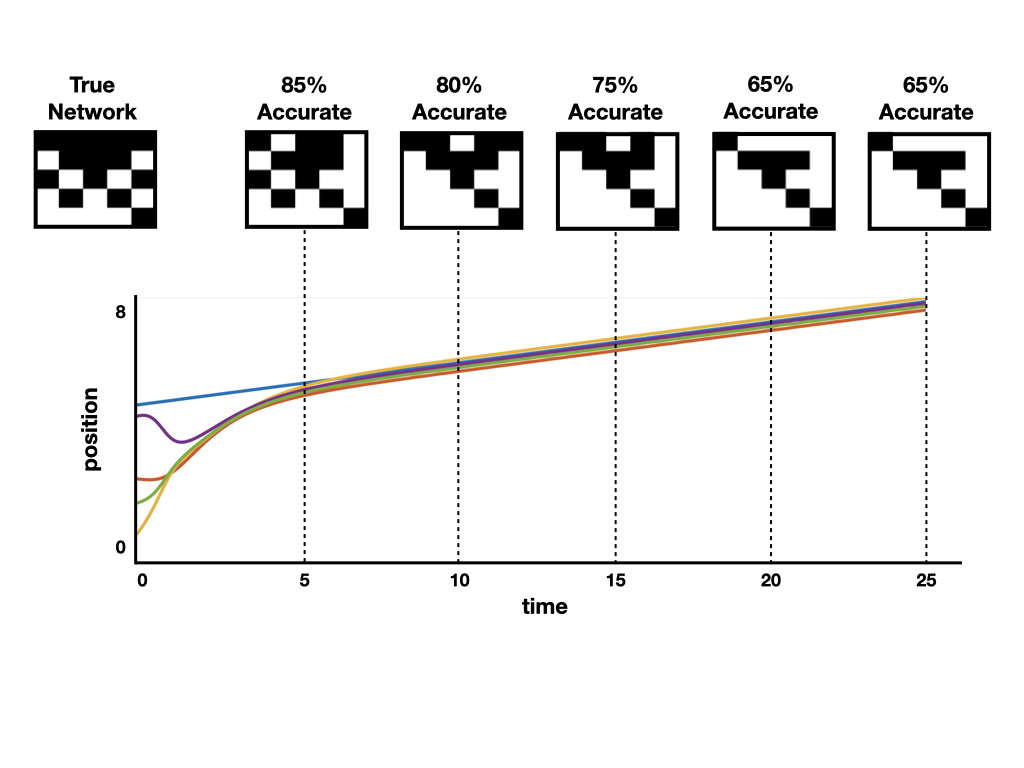}%
\vspace*{-.5in}
\caption{GC network prediction of a 5 node Kuramoto oscillator network with a coupling strength of 5 when length of time series data analyzed is stopped at 5, 10, 15, 20 or 25 seconds. Each of the colored lines is a trajectory of one of the five oscillators. Note that the network inferences are best when the length of the simulation passed into GC is close to the transient time length of 5 seconds.}
\label{fig:kuramotoGCTransients}
\end{figure}

In the case of the mass-spring system, the damped harmonic oscillators eventually tend to their equilibrium positions and in the case of the Kuramoto model, the oscillators synchronize and rotate with the same angular velocity. After energy imparted on the system is lost, the dynamics of the coupled oscillators become uninformative for the purposes of network reconstruction as shown in Fig. \ref{fig:kuramotoGCTransients}. The perturbation inference approach discussed here takes advantage of the transient time dynamics of a system directly after it is perturbed, while disregarding the complex nonlinear interactions that can arise later in time.

\begin{figure}
\includegraphics[width=8cm]{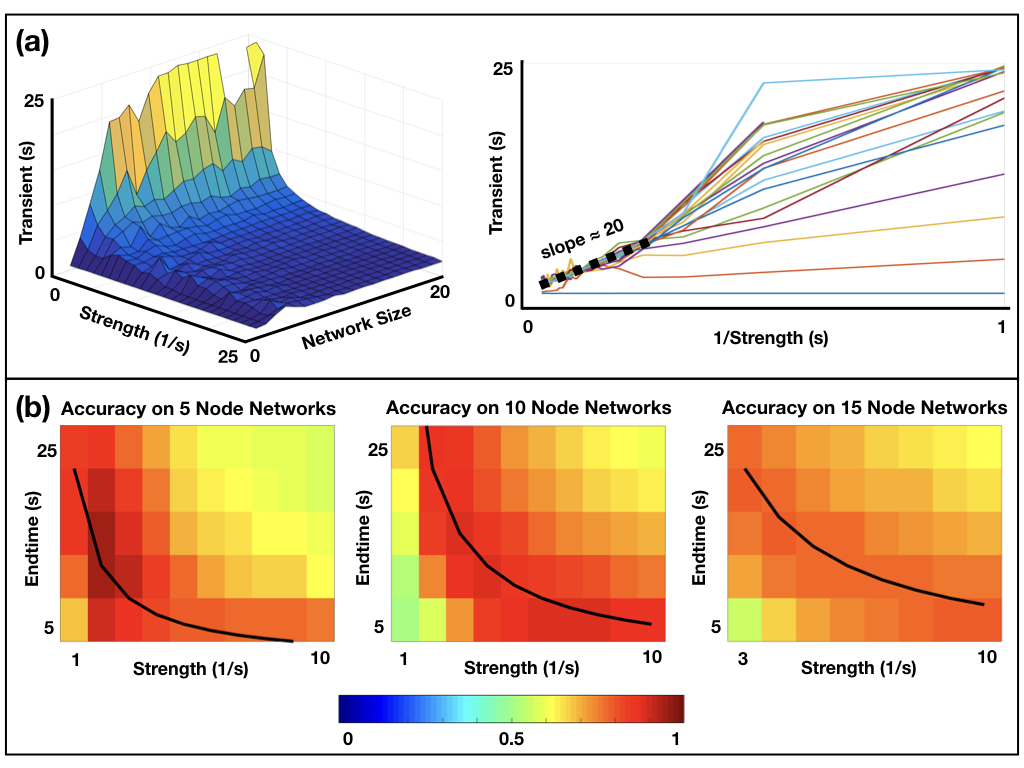}%
\caption{\textbf{(a)} Plot of transient time length of Kuramoto network vs coupling strength and network size. The transient time length of the system depends inversely on the coupling strength, and this dependence stays approximately the same for all networks larger than 5 nodes. Right graph plots the transient time length vs the inverse coupling strength of Kuramoto networks ranging from 1 to 20 nodes in size. A cluster of linear regions can be seen in the graph with a slope of approximately 20 and this cluster corresponds to all networks 5 nodes or larger. Therefore, for moderately large networks the transient time length is approximately 20 over the coupling strength of the system.
\textbf{(b)} GC accuracy over Kuramoto networks with sizes $n = 5, 10, 15$ and varying simulation endtimes and coupling strengths. Black line in heatmap plots corresponds to $T(n, K) = \frac{4.5n}{K}$ (optimal simulation time length is equal to 4.5 times the network size over the coupling strength). Note that for all network sizes, GC achieves a maximum accuracy around $T(n, K)$.}
\label{fig:kuramotoGCTransientsStrength}
\end{figure}

Fig.~\ref{fig:kuramotoGCTransients} shows an example of a networked, nonlinear dynamical system of Kuramoto oscillators.  We have already shown in Fig.~\ref{fig:inferenceAccuracy} that GC fails when simply presented with a given set of time-series data. But by considering the transient dynamics, and in fact, only using the first 5 seconds of transient dynamics, the GC accuracy is 85\% in reconstructing the true network. This shows that the transients contain most of the information required for inferring causality. By sampling for longer times, the dynamics fall into an attractor which obfuscates the true network architecture.  A more careful analysis of the Kuramato oscillator system is considered in Fig.~\ref{fig:kuramotoGCTransientsStrength}. The analysis of this plots shows the dependency of the sampling time on the coupling strength, which is directly related to the transient time exhibited in the dynamics.  Thus the time sampled of the dynamics, and its transient behavior, is critical in determining an accurate representation of causal relations. A  more detailed assessment of GC when the oscillator system is physically perturbed (as opposed to restarted with new initial conditions) is done in Sec.~IV.

\newpage

\subsection{Perturbation Cascade Inference}

\begin{figure}
    \includegraphics[width=9cm]{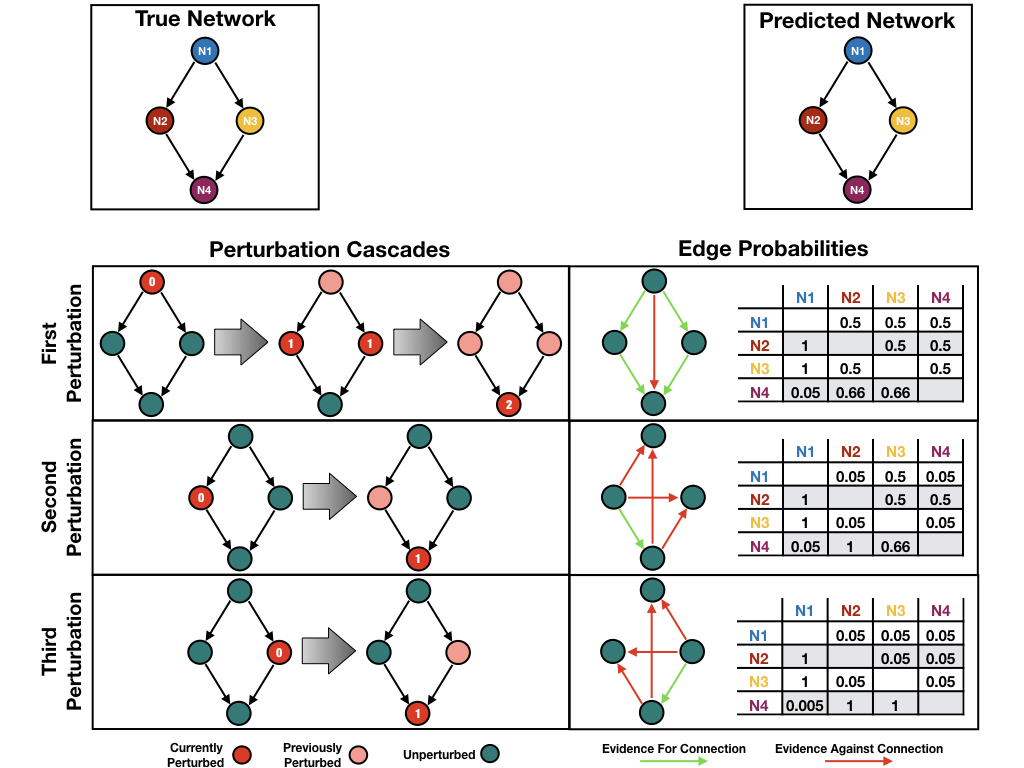}
    \caption{A simple example illustrates the strategy to infer the structure of a network of four coupled nodes using systematic perturbations. We make the important assumption that when a node is activated, it perturbs only its nearest neighbors at the next point in time. In the first perturbation, N1 is activated which causes N2 and N3 to become perturbed which in turn cause N4 to become perturbed. Originally, all edges in the network are given a probability of 0.5. This cascade of perturbations provides new evidence that N1 causes N2 and N3, either N2 or N3 causes N4, and N1 does not cause N4. Therefore, the edges from N1 to N2 and N3 are updated to 1, the edges from N2 and N3 to N4 are updated to 0.66, and the edge from N1 to N4 is updated to 0.05. This process is repeated for the second perturbation of N2 and the third perturbation of N3 with the adjacency matrix of edge probabilities updated each time. The final probability matrix is thresholded at 0.5 which results in the fully accurate network prediction shown in the top right.}
    \label{fig:inferenceWithPerts}
\end{figure}

The PCI method developed here physically forces nodes in the network to infer its underlying structure. It relies on the observation that different nodes are activated in time as information and perturbations spread across a dynamical network.  PCI learns the length of the shortest path (e.g. distance) from every node in the network to the initially perturbed node.  Indeed, PCI learns the distances of every node in the network from each perturbed node, and applies this information to reconstruct the underlying graph connectivity. To implement this approach on the mass-spring and Kuramoto oscillator networks studied in this paper, we apply impulse forcings at particular network nodes and study the transient time dynamics of the oscillators in the short time window after the perturbation was applied. We correlate the oscillator trajectories to sort all the oscillators in order of their activations in time and use these node orderings to determine the probabilities of all edges in the network. It is important to stress that this methodology is only applicable when the coupling strengths between all pairs of neighboring nodes in the oscillator network are approximately equal. Figure \ref{fig:inferenceWithPerts} shows a diagram of this method on a four node network.

To reverse-engineer the structure of a network, we must understand how to probabilistically infer its edges from information about the distance of every node from a set of perturbed nodes. Let's assume that after perturbing node $p$ in the network $\mathbf{A}$ we infer the length of the shortest directed path from $p$ to any node $v \in [n]$, denoted by $\dist(p, v)$. By convention, if there is no directed path from $p$ to $v$ then we write $\dist(p, v) = \infty$. Now we define the set of all nodes which have a directed shortest path of length $k$ from node $p$
\begin{equation}
    D_k(p) = \{v: \dist(p, v) = k\}.
\end{equation}
By this definition, $D_1(p)$ is the set of all the nearest neighbors of $p$ (i.e. nodes linked by a directed edge from $p$) and $D_\infty(p)$ contains all the nodes which are disconnected from $p$. Notice that any node $y \in D_{k+1}(p)$ must be connected to some node $x \in D_k(p)$ by an outgoing edge from $x$ (i.e. caused by $x$).  Furthermore, any node $x \in D_\infty(p)$ cannot be caused by/linked to any node in $D_k(p)$ where $k < \infty$. Therefore, if by perturbing node $p$ we can learn $\dist(p, v)$ for every node $v$ in the network, then we discover that:
\begin{enumerate}
    \item All nodes caused by $p$ lie in $D_1(p)$.\\
    \item Every node in $D_{k+1}(p)$ is caused by a node in $D_k(p)$.\\
    \item Every node in $D_k(p)$ for $k < \infty$ does not cause any of the nodes in $D_\infty(p).$
\end{enumerate}

As we perturb more nodes, we would like to assign a probability to each edge in the network that it actually exists in the underlying graph. We start with an uninformative prior where every edge in the network $x \to y$ has a $\mathbb{P}(x \to y) = 0.5$ probability of existing. Then for every subsequent perturbation, we update the network edge probabilities as follows:

Suppose we have performed a perturbation of node $p$ in the network and have learned the sets $D_1(p), D_2(p), ..., D_\infty(p)$ by sorting all of the oscillators in order of their activation times. Now we take any nodes $x, y$ from our network and update the probability that the edge $x \to y$ exists as follows:\\

\begin{algorithm}[t]
\SetAlgoLined
\KwResult{Returns edge probability matrix A}
 \tcp{A(x, y) probability x causes y}
 A = ones(n, n)/2\;
 \BlankLine
 perturbNodes = [1, 2, ...]\;
 \For{$p \in \mathrm{perturbNodes}$}{
  perturb node $p$\;
  observe sets $D_1(p), ..., D_\infty(p)$\;
  \For{$x \in [n]$}{
   \tcp{Case 1}
   \If{$x \in D_\infty(p)$}{
    continue\;
   }
   \BlankLine
   \tcp{Case 2}
   $k = \dist(p, x)$\;
   \For{$y \notin \bigcup_{m=1}^{k+1} D_m(p) \cup D_\infty(p)$}{
    $A(x, y) \leftarrow A(x, y) / 10$\;
   }
   \BlankLine
   \tcp{Case 3}
   \For{$y \in D_{k+1}(p)$}{
    $c = \prod_{v \in D_k(p)} (1 - A(v, y))$\;
    $A(x, y) \leftarrow A(x, y) / (1 - c)$\;
   }
  }
 }
 \caption{Perturbation Cascade Inference \label{alg:one}}
\end{algorithm}

\textbf{Case 1:}
If $x \in D_\infty(p)$ then $x$ is not caused by $p$ and we have not learned any information about the edge $x \to y$ so $\mathbb{P}(x \to y)$ is not updated.\\

\textbf{Case 2:} If $x \in D_k(p)$ for some $0 \leq k < \infty$ and $y \in D_m(p)$ where $k + 1 < m \leq \infty$ then we know from the discussion above that $x$ does not cause $y$ so we penalize the prior probability $\mathbb{P}(x \to y)$ by dividing it by 10.\\

\textbf{Case 3:} If $x \in D_k(p)$ and $y \in D_{k + 1}(p)$, where $0 \leq k < \infty$, we apply Bayes's rule to update the probability $\mathbb{P}(x \to y)$. Since $y \in D_{k + 1}(p)$, we know that $y$ is caused by at least one of the nodes in $D_k(p)$ (not necessarily by $x$). Each node $v \in D_k(p)$ has \textit{a prior} probability $\mathbb{P}(v \to y)$ that it causes $y$. Using Bayes's rule, we know that the conditional probability that the edge $x \to y$ exists given the information that at least one of the edges $v \to y$ exists for $v \in D_k(p)$ is
\begin{eqnarray}\label{eq:13}
&&\mathbb{P}\Big(\!x \to y \ \!\Big|\!\ \exists v \!\in\! D_k(p) \!\text{ s.t. }\! v \to y \!\Big) \nonumber \\
&&\hspace*{.3in} = \! \frac{\mathbb{P}\Big(\exists v \!\in\! D_k(p) \!\text{ s.t. }\! v \to y \ \!\Big|\!\ x \to y\Big)\mathbb{P}(x \to y)}{\mathbb{P}\Big(\!\exists v \!\in\! D_k(p) \text{ s.t. } v \to y\!\Big)} \nonumber \\
&& \hspace*{.3in} =\! \frac{\mathbb{P}(x \to y)}{1 - \prod_{v \in D_k(p)}\! \Big(\!1 \!-\! \mathbb{P}(v \to y)\!\Big)}.
\end{eqnarray}
These computations produce a probability of a causal connection between one node and another. The entire algorithm is summarized in Algorithm~\ref{alg:one}.

In the discussion above, we did not mention how to infer the distance sets $D_k(p)$ after perturbing node $p$ in the network. The construction of these sets is an entirely data-driven problem that heavily depends on the time-dynamics of the system being studied. In the case of the simulated mass-spring and Kuramoto oscillator networks, we use correlations and variances between oscillator trajectories respectively to estimate which nodes lie in each distance set $D_k(p)$ by studying at what times they become activated.
The effectiveness of PCI is illustrated in the following section.

\section{Computational Results:  Coupled Oscillators}
\label{sec:results}

\begin{figure*}
\includegraphics[width=15cm]{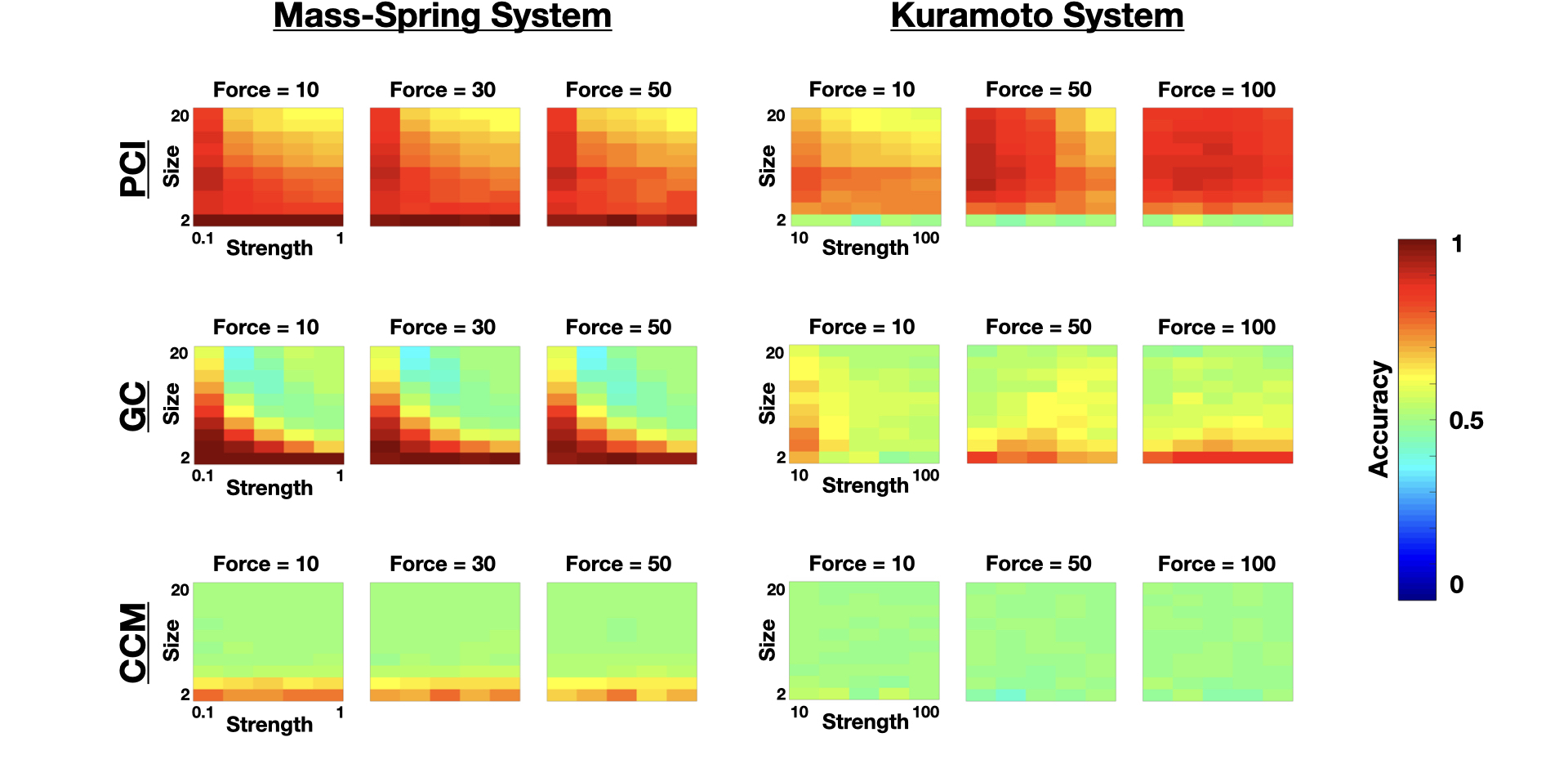}%
\caption{Accuracy of PCI, GC, and CCM on perturbed mass-spring and Kuramoto oscillator ER networks (0.5 connection probability) with varying sizes, perturbation forcing magnitudes, and coupling strengths. For the mass-spring system, PCI achieves $\approx$80\% accuracy for all networks of at most 10 nodes and this accuracy increases slightly if the force magnitude by which the nodes are perturbed is increased from 10N to 50N. For larger network sizes, PCI requires lower coupling strengths in order to reach the same predictive accuracy. As an example, for the 20 node Erdos-Renyi networks, PCI achieves this accuracy only when all springs in the mass-spring network have spring constants equal to 0.1. For the Kuramoto oscillator system PCI requires large coupling and forcing strengths $(>20)$ in order to correctly infer networks of 2 nodes or larger. Therefore, if the coupling and forcing of the system are large, then PCI successfully infers Kuramoto oscillator networks with an accuracy of 80-90\%. GC is capable of inferring network structure of a 10 node mass-spring system with 85\% accuracy when the coupling strengths are 0.1 or smaller. Similar to PCI, its performance drops below 80\% for larger coupling strengths and network sizes and is not significantly affected by the magnitude with which the nodes are perturbed. GC consistently gets below 75\% accuracy on Kuramoto oscillator networks with 3 nodes or larger. Note that even with these perturbed systems, CCM is incapable of accurately predicting any structure and for networks of larger than 3 nodes, achieves an accuracy of close to 50\%.}
\label{fig:pertSFS}
\end{figure*}

In this section, we evalute the performance of GC, CCM, and PCI on mass-spring and the nonlinear Kuratomo oscillator systems where the oscillators are \textit{physically forced}.
Fig.~\ref{fig:pertSFS} gives an extensive comparison of the models as a function of the network size, coupling strength between nodes and number of forcings (perturbations). The networks are randomly generated Erdos-Renyi (ER) graphs with a 0.5 probability of connection. Fig.~\ref{fig:pertSFS} shows that CCM fails to accurately infer causal connections even when perturbations are allowed, much like what has been found in \cite{monster2017,monster2018}. Perturbed GC has a range of parameter space where accurate causal inference can be established, although its validity is especially limited for the nonlinear oscillators. PCI has strong potential for accurate causal inference, especially as the number of perturbations to the system is increased.

\begin{figure*}
\includegraphics[width=15cm]{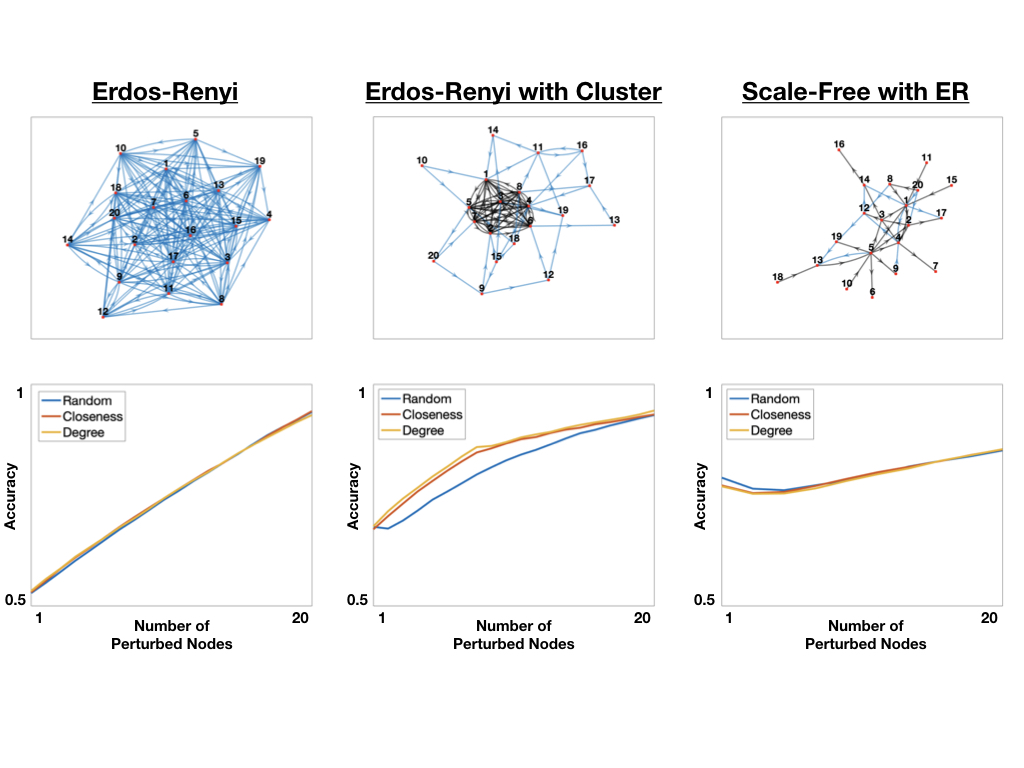}%
\vspace*{-.5in}
\caption{The accuracy of PCI on inferring 20 node mass-spring networks, when nodes are perturbed randomly or in order of their outdegree or centrality. We perform this experiment on three different network types: Erdos-Renyi, Erdos-Renyi with cluster, and scale-free. We observe one case (namely on the clustered ER graphs), where perturbing less than half of the network nodes (8/20) results in an optimistic 85\% inference accuracy.}
\label{fig:PCIAccVSPerts}
\end{figure*}

Given the success of the PCI method, we explore its use on larger 20 node mass-spring networks as shown in Fig.~\ref{fig:PCIAccVSPerts}. In these experiments, we infer the structure of an Erdos-Renyi network with 0.5 connection probability, a 0.5 probability Erdos-Renyi network with a fully-connected cluster/clique of 8 nodes, and a scale-free network built by the Barab\'{a}si-Albert model~\cite{albert2002statistical}. For each graph, we try three different orders of node perturbations: random order, by decreasing order of outdegree, and by decreasing order of outcloseness. In a real experiment, we would never have a full understanding of the degree or closeness statistics of all nodes in the network. However, this experiment is meant to simulate the ``best" possible scenario for inference if we have prior knowledge about the importance or centrality of nodes in the graph.
We do not consider either perturbed GC or CCM as their performance does not scale to 20 node networks and does not significantly depend on the order of node perturbations.

PCI is shown to be a viable technique for extracting accurate causal relations provided enough perturbations of the system are allowed. The causal relations are nicely extracted for ER networks, and with ER networks containing clusters and scale-free structure. In particular, on the ER cluster network we observe that if we perturb in decreasing order of degree or closeness instead of randomly we can reach 85\% prediction accuracy with only 8 out of the 20 nodes perturbed. This shows that in some instances, PCI can infer certain network structures with a small set of properly placed perturbations, similar to the way that messages efficiently spread throughout a network when they are transmitted from highly-connected hubs.

\section{Conclusion}

In this paper, we showed that perturbing components of the network, or by observing many unique transient dynamics, a large number of potential networks can be disambiguated, achieving reconstruction of unique and accurate network structures. Specifically, we demonstrated that our PCI method, along with GC with perturbations and transients, is capable of inferring linear (mass-spring) and nonlinear (Kuramoto oscillator) networked dynamical systems. 
Our proposed PCI method demonstrated consistently strong performance in inferring causal relations for small (2--5 node) and large (10--20 node) networks for both linear and nonlinear systems. Perturbed GC is capable of inferring smaller networks under low coupling strength regimes, while methods such as CCM did not infer the structure of any oscillator networks. Our analysis suggests that the ability to apply a large and diverse set of perturbations/actuations to the network, in either a targeted or random way, is critical for successfully and accurately determining causal relations and disambiguating between various viable networks.

Beyond data analysis and modeling, the problem of network inference asks a foundational question: How can a complex system be described by a small, interpretable set of causal relationships among its components?
In many scientific and engineering applications, studying transient dynamics and impulse responses have long provided insight on the desired causal structure. 
We note, however, that to obtain reasonable network inferences and interpret them, much care must be taken to choose and observe only the important nodes of multi-component networks rather than analyzing the entire system in all of its complexity.  
Importantly, there remain many open fundamental mathematical questions for future study.
These include how to optimally place activations and observations of network nodes for inference, what types of network dynamics contain unique information about their connectivity, and above all, whether it is possible to formulate a consistent and unifying theory of causality for time dependent systems.

\section*{Acknowledgements}
GS would like to acknowledge support from the Mary Gates Research Foundation and the National Science Foundation Graduate Research Fellowship under Grant No. 174530. 
JNK acknowledges support from the Air Force Office of Scientific Research (AFOSR) grant FA9550-17-1-0329.
BWB acknowledges support from the Washington Research Foundation and Air Force Research Lab (AFRL) grant FA8651-16-1-0003.

All code and additional experiments can be found at \href{https://github.com/sgstepaniants/netinf}{https://github.com/sgstepaniants/netinf}.

\section*{Appendix}

The dynamical models used to evaluate causal relations are of two type:  linear and nonlinear.  The linear model is a standard mass-spring system where masses interact through springs and Hooke's law.  The nonlinear system is a set of nonlinear oscillators known as the Kuramoto system. Details are given in each section below.

\subsection{Mass-Spring System}
Mass-spring networks have numerous applications in a variety of disciplines including modeling of deformable objects in computer graphics \cite{nealen2005}, molecular dynamics in complex polymer materials \cite{zhang2014} and organ simulations for surgical procedures \cite{meier2005}. The system of ODEs that governs the movement of $n$-coupled oscillators in a directed mass-spring network is
\begin{equation}\label{eq:5}
m\frac{d^2\mathbf{x}}{dt^2} = k\mathbf{M}\mathbf{x} - c\frac{d\mathbf{x}}{dt} + \mathbf{f}(t),
\end{equation}
where we assume the mass $m$ and damping constant $c$ of every oscillator is the same and the spring constant $k$ is also the same for every spring in the network. $\mathbf{f}(t)$ is the forcing function for every oscillator. The binary matrix $\mathbf{M}$ defines the connectivity of the network. We impose fixed boundary conditions such that the first and last oscillators are attached to fixed walls by springs. For all other oscillators in the network, they cannot be directly connected to the boundary walls. Therefore, $\mathbf{M}$ has dimension $n + 1 \times n + 1$ and has the form
\begin{equation}
\mathbf{M} =
\begin{bmatrix}
    0 & 1 & 0 & \hdots & 0\\ \cline{2-4}
    1 & \multicolumn{3}{|c|}{\multirow{3}{*}{\scalebox{2}{$\mathbf{A}$}}} & \vdots\\
    0 & \multicolumn{3}{|c|}{} & 0\\
    \vdots & \multicolumn{3}{|c|}{} & 1\\ \cline{2-4}
    0 & \hdots & 0 & 1 & 0
\end{bmatrix},
\end{equation}
where the $n \times n$ binary adjacency matrix $\mathbf{A}$ represents the connectivity of all of the nodes in the network with the exception of the fixed boundary walls whose connectivity is predetermined. If $A_{ji} = 1$ then oscillator $i$ is connected to oscillator $j$ by a `directed spring' (i.e. oscillator $i$ can force oscillator $j$) and if $A_{ji} = 0$ then oscillator $i$ has no direct influence on oscillator $j$. Since causality is a directed relationship (i.e. $X$ causes $Y$ but $Y$ might not cause $X$), we do note impose restrictions that $\mathbf{A}$ must be symmetric.

\subsection{Kuramoto System}
The Kuramoto model proposed by Yoshiki Kuramoto is one of the most well-studied systems of nonlinear coupled oscillators \cite{kuramoto1975, kuramoto2003, strogatz2000kuramoto}. It is a canonical system for studying quasiperiodic dynamics, synchronization, and chaos and has found practical applications in a variety of areas in physics \cite{strogatz1998, strogatz1996, kourtchatov, jiang}, biology \cite{liu1997, buck1988, walker1969}, and medical sciences \cite{peskin1975, michaels1987}. Equation (\ref{eq:7}) describes the dynamics of $n$ coupled Kuramoto oscillators with a forcing term.
\begin{equation}\label{eq:7}
\frac{d\theta_i}{dt} = \omega_i + \frac{K}{n}\sum_{j=1}^n A_{ij}\sin(\theta_j - \theta_i) + f_i(t), \quad 1 \leq i \leq n.
\end{equation}
The dynamics of the $i$th oscillator is described by its angle $\theta_i$ which has a natural frequency $\omega_i$. Each $A_{ij}$ is an entry of the $n \times n$ binary adjacency matrix $\mathbf{A}$ which represents the connectivity of all $n$ oscillators in the network. All oscillators in the network are coupled to their adjacent neighbors by the same coupling strength $K$. For any Kuramoto network there exists a bifurcation value $K_c$ where for all $0 \leq K < K_c$ the oscillator trajectories are unsynchronized and for $K > K_c$ clusters of oscillators synchronize and eventually phase lock. This bifurcation value depends on various properties including the topology of the network and the distribution from which the natural frequencies $\omega_i$ are sampled from. Also, the Kuramoto model exhibits chaotic dynamics for networks of four or more nodes \cite{maistrenko2005}. Therefore, in all our experiments we test a sufficiently wide range of coupling strengths and network sizes to analyze how network inferences vary in synchronized, unsynchronized, chaotic, and non-chaotic regimes.

\bibliography{references.bib}

\end{document}